\documentclass[10pt]{article}
\usepackage{amsfonts}
\usepackage{amsmath,amssymb,amscd,verbatim,amsthm}
\usepackage{amsmath}
\usepackage{amssymb}

\begin{document}

\begin{center}
{\Large \textbf{A short note on the paper \textquotedblleft Remarks on
Caristi's fixed point theorem and Kirk's problem\textquotedblright}}\\[0.3in%
]
\textbf{Wei-Shih Du}\bigskip 

{\large {\small Department of Mathematics, National Kaohsiung Normal
University, Kaohsiung 802, Taiwan}\bigskip }

{\large {\small E-mail address: wsdu@nknucc.nknu.edu.tw}\\[0.3in]
}
\end{center}

\vspace{0.3cm} \hrule\vspace{0.2cm}\bigskip

{\small \noindent \textbf{Abstract:} In this paper, we demonstrate that Li's
fixed point theorems are indeed equivalent with the primitive Caristi's
fixed point theorem, Jachymski's fixed point theorems, Feng and Liu's fixed
point theorems, Khamsi's fixed point theorems and others.}\\[0.05in]


\\[0.005in]
\vspace{0.2cm} \hrule\vspace{0.3cm}\bigskip \bigskip

Let $(X,d)$ be a metric space. Recall that an operator $T:X\rightarrow X$ is
said to be a \textit{Caristi type} mapping [1] provided that the following
condition is satisfied
\begin{equation*}
\eta \left( d(x,Tx)\right) \leq \varphi (x)-\varphi (Tx)\text{, }\forall
\text{ }x\in X\text{,}
\end{equation*}%
where $\eta :[0,+\infty )\rightarrow
\mathbb{R}
$ and $\varphi :X\rightarrow (-\infty ,+\infty ]$ are functions. \bigskip

In 1976, Caristi [2] proved the following famous fixed point theorem:
\bigskip

\noindent \textbf{Theorem CA. [2]}\quad Let $(X,d)$ be a complete metric
space and $\varphi :X\rightarrow
\mathbb{R}
$ be a lower semicontinuous and bounded below function. Suppose that $%
T:X\rightarrow X$ is a single-valued mapping satisfying
\begin{equation*}
d(x,Tx)\leq \varphi (x)-\varphi (Tx)\text{, }\forall \text{ }x\in X\text{,}
\end{equation*}%
that is $T$ is a Caristi type mapping with $\eta (t)=t$, $t\in \lbrack
0,\infty )$. Then $T$ has a fixed point in $X$.\bigskip

It is well known that the primitive Caristi's fixed point theorem is
equivalent to Ekeland's variational principle [3] and to Takahashi's
nonconvex minimization theorem [4, 5]. A number of generalizations in
various different directions of these results in metric (or quasi-metric)
spaces and more general in topological vector spaces have been investigated
by several authors in the past; see [6-19, 22-25] and references therein.
But it is worth to mention that some generalizations were real logical
equivalent with the original theorems; see e.g. [19, 23-25].

A problem rasied by Kirk [20-23, 29, 30] asked whether a Caristi type
mapping $T$ for a suitable function $\eta $ has a fixed point. In fact the
original Kirk's question was stated when $\eta (t)=t^{p}$ for some $p>1$. Some negative answers to this problem were 
given; for more detail, one can see [21-23, 25]. Very recently, Li [1]
investigated the existence fixed points for Caristi type mappings which
partially answered Kirk's problem and improved Caristi's fixed point
theorem, Jachymski's fixed point theorems [9, 12], Feng and Liu's fixed
point theorems [14], Khamsi's fixed point theorems [23] and others.

The following new fixed point theorems are the main results in [1].\bigskip

\noindent \textbf{Theorem 1 [1].}\quad Let $(X,d)$ be a complete metric
space. Suppose that $\eta :[0,+\infty )\rightarrow
\mathbb{R}
$ with $\eta (0)=0$, $\varphi :X\rightarrow
\mathbb{R}
$ is lower semicontinuous on $X$, and there exist $x_{0}\in X$ and two real
numbers $a$, $\beta \in
\mathbb{R}
$ such that
\begin{equation*}
\varphi (x)\geq ad(x,x_{0})+\beta \text{.}
\end{equation*}%
If one of the following conditons is satisfied:

\begin{enumerate}
\item[(i)] $a\geq 0$, $\eta $ is nonnegative and nondecreasing on $%
W=\{d(x,y):x,y\in X\}$, and there exist $c>0$ and $\varepsilon >0$ such that
\begin{equation*}
\eta (t)\geq ct\text{, }\forall \text{ }t\in \{t\geq 0:\eta (t)\leq
\varepsilon \}\cap W\text{;}
\end{equation*}

\item[(ii)] $a<0$, $\eta (t)+at$ is nonnegative and nondecreasing on $W$,
and there exist $c>0$ and $\varepsilon >0$ such that
\begin{equation*}
\eta (t)+at\geq ct\text{, }\forall \text{ }t\in \{t\geq 0:\eta (t)+at\leq
\varepsilon \}\cap W\text{.}
\end{equation*}
\end{enumerate}

\noindent Then each Caristi type mapping $T:X\rightarrow X$ have a fixed
point in $X$.\bigskip

\noindent \textbf{Theorem 2 [1].}\quad Let $(X,d)$ be a complete metric
space. Suppose that $\eta :[0,+\infty )\rightarrow \lbrack 0,+\infty )$ with
$\eta (0)=0$, $\varphi :X\rightarrow
\mathbb{R}
$ is lower semicontinuous on $X$ and bounded below on each bounded subset of
$X$, and there exist $x_{0}\in X$ and a real number $a\in
\mathbb{R}
$ such that
\begin{equation*}
\liminf_{d(x,x_{0})\rightarrow +\infty }\frac{\varphi (x)}{d(x,x_{0})}>a%
\text{.}
\end{equation*}%
If one of the following conditons is satisfied:

\begin{enumerate}
\item[(i)] $a\geq 0$, $\eta $ is nonnegative and nondecreasing on $%
[0,+\infty )$, and
\begin{equation*}
\liminf_{t\rightarrow 0+}\frac{\eta (t)}{t}>0\text{;}
\end{equation*}

\item[(ii)] $a<0$, $\eta (t)+at$ is nonnegative and nondecreasing on $%
[0,+\infty )$, and
\begin{equation*}
\liminf_{t\rightarrow 0+}\frac{\eta (t)}{t}\geq -a\text{;}
\end{equation*}
\end{enumerate}

\noindent Then each Caristi type mapping $T:X\rightarrow X$ have a fixed
point in $X$.\bigskip

In [1], Li had shown the following.\bigskip

\begin{enumerate}
\item[$\bullet $] Theorem CA $\Rightarrow $ Theorem 1.

\item[$\bullet $] Theorem 1 $\Rightarrow $ Theorem 2.
\end{enumerate}

\bigskip

It is obvious that each of Theorems 1 and 2 implies Theorem CA (the
primitive Caristi's fixed point theorem), so we can obtain the following
very important result.\bigskip

\noindent \textbf{Theorem D.}\quad Theorem CA, Theorem 1 and Theorem 2 are
equivalent.\bigskip

\noindent \textbf{Remark.}\quad

\begin{enumerate}
\item[(a)] Li also gave [1, Corollary 1] by using Theorem 2 and he point out
that some known extensions of Caristi's fixed point theorem established by
Downing and Kirk [6, 7], Jachymski [9, 12], Feng and Liu [14], and Khamsi
[23] are special cases of Theorem 2 (for more detail, see [1, Remark 2] and
[1, Remark 3]). So, by Theorem D, we know that they are real logical
equivalent.

\item[(b)] By a similar argument as Li's results, we can obtain easily that
Li's fixed point theorems for other weak distances ($w$-distances [5, 8, 16]
or $\tau $-distances [10, 13] or $\tau $-functions [15, 17, 31] and so on)
are equivalent with some weak distance variants of Caristi's fixed point
theorem, Ekeland's variational principle and Takahashi's nonconvex
minimization; see, e.g., [5, 8, 10, 15-17, 19].
\end{enumerate}

\bigskip \bigskip \bigskip

\noindent {\Large \textbf{References}}

\begin{enumerate}
\item[{[1]}] Z. Li, Remark on Caristi's fixed point theorem and Kirk's
problem, Nonlinear Analysis (2010), doi:10.1016/j.na.2010.07.048.

\item[{[2]}] J. Caristi, Fixed point theorems for mappings satisfying
inwardness conditions, Trans. Amer. Math. Soc. 215 (1976), 241-251.

\item[{[3]}] I. Ekeland, Remarques sur les probl\'{e}mes variationnels, I, C.
R. Acad. Sci. Paris S\'{e}r. A-B. 275 (1972), 1057-1059.

\item[{[4]}] W. Takahashi, Existence theorems generalizing fixed point
theorems for multivalued mappings, in: M.A. Th\'{e}ra, J.B. Baillon (Eds),
Fixed point theory and Applications, Pitmam Research Notes in Mathematics
Series, Vol. 252, Longmam Sci. Tech., Harlow, 1991, pp. 397-406.

\item[{[5]}] W. Takahashi, Nonlinear Functional Analysis: Fixed point theory
and its applications, Yokohama Publishers, Yokohama, Japan, 2000.

\item[{[6]}] D. Downing, W.A. Kirk, A generalization of Caristi's theorem
with applications to nonlinear mapping theory, Pacific J. Math. 69 (1977)
339-345.

\item[{[7]}] D. Downing, W.A. Kirk, Fixed point theorems for set-valued
mappings in metric and Banach spaces, Math. Japon. 22 (1977) 99-112.

\item[{[8]}] O. Kada, T. Suzuki, W. Takahashi, Nonconvex minimization
theorems and fixed point theorems in complete metric spaces, Math. Japon. 44
(1996), 381-391.

\item[{[9]}] J.R. Jachymski, Caristi's fixed point theorem and selection of
set-valued contractions, J. Math. Anal. Appl. 227 (1998), 55-67.

\item[{[10]}] T. Suzuki, Generalized distance and existence theorems in
complete metric spaces, J. Math. Anal. Appl. 253 (2001), 440-458.

\item[{[11]}] J.S. Bae, Fixed point theorems for weakly contractive
multivalued maps, J. Math. Anal. Appl. 284 (2003), 690-697.

\item[{[12]}] J.R. Jachymski, Converses to fixed point theorems of Zeremlo
and Caristi, Nonlinear Anal. 52 (2003) 1455-1463.

\item[{[13]}] T. Suzuki, Generalized Caristi's fixed point theorems by Bae
and others, J. Math. Anal. Appl. 302 (2005), 502-508.

\item[{[14]}] Y.Q. Feng, S.Y. Liu, Fixed point theorems for multi-valued
contractive mappings and multi-valued Caristi type mapping, J. Math. Anal.
Appl. 317 (2006) 103-112.

\item[{[15]}] L.-J. Lin, W.-S. Du, Ekeland's variational principle, minimax
theorems and existence of nonconvex equilibria in complete metric spaces, J.
Math. Anal. Appl. 323 (2006) 360-370.

\item[{[16]}] L.-J. Lin, W.-S. Du, Some equivalent formulations of
generalized Ekeland's variational principle and their applications,
Nonlinear Anal. 67 (2007) 187-199.

\item[{[17]}] L.-J. Lin, W.-S. Du, On maximal element theorems, variants of
Ekeland's variational principle and their applications, Nonlinear Anal. 68
(2008) 1246-1262.

\item[{[18]}] L.-J. Lin, W.-S. Du, Systems of equilibrium problems with
applications to new variants of Ekeland's variational principle, fixed point
theorems and parametric optimization problems, J. Global Optim. 40 (2008),
663--677.

\item[{[19]}] W.-S. Du, On some nonlinear problems induced by an abstract
maximal element principle, J. Math. Anal. Appl. 347 (2008) 391-399.

\item[{[20]}] J. Caristi, Fixed point theory and inwardness conditions,
Applied Nonlinear Analysis, Academic Press (1979), 479-483.

\item[{[21]}] Nguyen Dinh Dan, Remark on a question of Kirk about Caristi's
fixed point theorem, Math. Operationsforsch. Statist. Ser. Optim. 12 (1981),
no. 2, 177-179.

\item[{[22]}] J. S. Bae and S. Park, Remarks on the Caristi-Kirk fixed point
theorem, Bull. Korean Math. Soc. 19 (1983), 57-60.

\item[{[23]}] M.A. Khamsi, Remarks on Caristi's fixed point theorem,
Nonlinear Anal. 71 (2009) 227-231.

\item[{[24]}] A. Latif, Generalized Caristi's fixed point theorems, Fixed
Point Theory and Applications, Volume 2009 (2009), Article ID 170140, 7 pages, doi:10.1155/2009/170140.

\item[{[25]}] R.P. Agarwal, M.A. Khamsi, Extension of Caristi's fixed point
theorem to vector valued metric spaces, Nonlinear Analysis (2010),
doi:10.1016/j.na.2010.08.025.

\item[{[26]}] T.Q. Bao, P.Q. Khanh, Are several recent generalizations of
Ekeland's variational principle more general than the original principle?
Acta Math. Vietnam. 28 (2003) 345--350.

\item[{[27]}] M. Turinici, Zhong's variational principle is equivalent with
Ekeland's, Fixed Point Theory 6 (1) (2005) 133-138.

\item[{[28]}] W.-S. Du, On Latif's fixed point theorems, Taiwanese Journal of
Mathematics, in press.

\item[{[29]}] W.A. Kirk, J. Caristi, Mapping theorems in metric and Banach
spaces, Bull. Acad. Polon. Sci. 23 (1975) 891-894.

\item[{[30]}] W.A. Kirk, Caristi's fixed-point theorem and metric convexity,
Colloq. Math. 36 (1976) 81-86.

\item[{[31]}] W.-S. Du, Some new results and generalizations in metric fixed
point theory, Nonlinear Anal. 73 (2010) 1439-1446.

\end{enumerate}

\end{document}